\documentclass{irmaart}
\usepackage[all]{xy}
\usepackage{amssymb}
\usepackage{amsmath}
\swapnumbers
\newtheorem{thm}[subsection]{Theorem}
\newtheorem{prop}[subsection]{Proposition}

\newtheorem{cor}[subsection]{Corollary}

\def\SS{\Sigma}

\def\dd{\delta}
\def\DD{\Delta}

\def\ss{\sigma}

\def\t{\otimes}

\def\ZZ{{\mathbb{Z}}}

\def\RR{{\mathbb{R}}}

\def\RR{{\mathbb{R}}}

\def\Hom{\mathrm{Hom}}

\def\Aut{\mathrm{Aut}}
\def\sgn{\mathrm{sgn}}

\def\Tot{\mathrm{Tot}}

\def\Hom{\mathrm{\rm Hom }}

\def\dim{\mathop{\rm dim }}
\def\id{\mathrm{ id }}

\def\KK{\mathbb{K}}

\def\LL{\mathcal{L}}
\def\SS{\mathcal{S}}

\def\DDb{{\bf \DD}}

\def\arbreA{\vcenter{\xymatrix@R=3pt@C=3pt{
&& \\
&*{}\ar@{-}[ur] \ar@{-}[ul] \ar@{-}[d]     &\\
&&
}}}

\def\Kzero{\xymatrix@R=4pt@C=4pt{
\\
\\
\\
{\bullet}\\
}}

\def\KunA{\xymatrix@R=4pt@C=4pt{
&&\\
&&\\
&&\\
*{}\ar@{-}[rr]&&*{}\\
}}

\def\TdeuxA{\xymatrix@R=4pt@C=4pt{
&&*{}\ar@{-}[dddll]\ar@{-}[dddrr]&&\\
&&&&\\
&&&&\\
*{}\ar@{-}[rrrr]&&&&*{}\\
}}

\def\TtroisA{\xymatrix@R=4pt@C=4pt{
&&&*{}\ar@{-}[ddddlll]\ar@{-}[ddddrrr]\ar@{-}[drrr]&&&\\
&&&&&&*{}\ar@{-}[ddd]\ar@{.}[dddllllll]\\
&&&&&*{}&\\
&&&&&&\\
*{}\ar@{-}[rrrrrr]&&&&&&*{}\\
}}

\begin{document}

\markboth{J.-L. Loday}{Free loop space and homology}
\address{Institut de Recherche Math\'ematique Avanc\'ee\\
    CNRS et Universit\'e de Strasbourg\\
    7 rue R. Descartes\\
    67084 Strasbourg Cedex, France\\
email:, \tt{loday@math.unistra.fr}
}

\title{Free loop space and homology}

\author{Jean-Louis Loday}

%\urladdr{www-irma.u-strasbg.fr/{$\sim$}loday/}

%\alttitle{} 
%\subjclass[2000]{ }
%\keywords{ }

%\date{\today}

\maketitle 

\begin{abstract} The aim of this paper is to explain the relationship between the (co)homology of the free loop space and the Hochschild homology of its singular cochain algebra. We introduce all the relevant technical tools, namely simplicial and cyclic objects, and we provide the various steps of the proofs, which are scattered around in the literature. This paper can be seen as a companion to the paper ``Cyclic homology and equivariant homology'' by J.D.S. Jones dealing with the cyclic case.
\end{abstract}
\vskip 1cm

\section*{Introduction} The purpose of this paper is to explain why Hochschild homology $HH$, and then cyclic homology $HC$, has anything to do with the notion of free loop space $\LL= \Hom(S^{1}, -)$. The answer, which is often hidden, or taken for granted, in the literature, is in the analysis of the simplest simplicial model of the circle $S^{1}$. This model, denoted by 
$S^{1}_{\cdot}$, has two non trivial cells: one in dimension $0$ (the base-point), one in dimension $1$. Then it turns out that $S^{1}_{n}$ is a set with $n+1$ elements. This is in accordance with the fact that the Hochschild complex of an associative algebra $A$ is $A\mapsto A^{\t n+1}$ in dimension $n$. Moreover the simplicial structure of $S^{1}_{\cdot}$ fits with the Hochschild boundary map $b$.

But there is more. Not only $S^{1}_{\cdot}$ has $n+1$ elements, but in fact it comes with a bijection $S^{1}_{n}\cong \ZZ/(n+1)\ZZ$ such that the group structure behaves well with the simplicial structure. All these properties are encompassed in the statement: $S^{1}_{\cdot}$ is a cyclic set (a discovery due to Alain Connes). On the homological side it gives rise to cyclic homology ; on the topological side it gives rise to the group structure of the circle and then to the $S^{1}$-structure of the free loop space.

Putting everything together one can construct isomorphisms:

$$ HH_{\bullet}(\SS^{* }X) \cong H^{\bullet}(\LL X)\quad \textrm{and} \quad HC_{\bullet}(\SS^{*} X) \cong H^{\bullet}_{S^{1}}(\LL X).$$
where $\SS^{* }X$ is the  algebra of singular cochains of the simply connected topological space $X$.
If $M$ is a simply connected manifold, then the de Rham algebra $\Omega^{*}M$ of $M$ is homotopy equivalent to $\SS^{* }M$, hence we get isomorphisms:

$$ HH_{\bullet}(\Omega^{*}M) \cong H^{\bullet}(\LL M)\quad \textrm{and} \quad HC_{\bullet}(\Omega^{*}M) \cong H^{\bullet}_{S^{1}}(\LL M),$$

Hochschild homology and cyclic homology are related by the ``periodicity exact sequence'' also known as ``Connes exact sequence''. It can be shown that, under the previous isomorphisms it corresponds exactly to the Gysin sequence associated to the action of the topological group $S^{1}=SO(2)$ on the free loop space.

There are many variations of these results, all well documented in the literature. Our aim in this paper, which does not contain anything new, is to provide an easy self-contained access to the aforementioned results.

\medskip

\noindent{\sc Contents}

1. Simplicial

2. Hochschild homology and free loop space

3. Cyclic

4. Cyclic homology and the free loop space

5. Conclusion and comments

\medskip

\noindent{\sc Acknowledgements.} Thanks to Luc Menichi and Jean-Claude Thomas for answering my questions. This paper would not exist without the constant and friendly determination of Alex Oancea.

%%%%%%%%%%%%%%%%%%%%%%%%%%%%%%%%%%%%%%%%%%%%%%%%%%
\section{Simplicial}\label{simplicial}
 In order to study manifolds up to homotopy, it is helpful to have a universal way of constructing them. Simplicialization is an efficient way since it permits us to apply the full strength of algebraic methods. In this section we recall the notions of (co)simplicial sets and (co)simplicial modules and we compare them with topological spaces and chain complexes respectively:
 
 $$\texttt{Top}\leftrightarrows \texttt{Simp Set} \to  \texttt{Simp Mod} \leftrightarrows \texttt{Chain Cplx}$$

\subsection{Simplices and the simplicial category $\DD$} The simplex of dimension $n$ is the convex hull of the points $(0, \ldots,0, \underbrace{1}_{i},0, \ldots, 0)\in \RR^{n+1}$ for $i$ ranging from $0$ to $n$. It is denoted by $\DDb^{n}$. For $n=0,1,2, 3$ we get respectively a point, an interval, a triangle, a tetrahedron: 

$$\Kzero \qquad \KunA \qquad \TdeuxA \qquad  \TtroisA $$

Opposite to the vertex number $i$ is the $i$th face which is a simplex of dimension $n-1$. The associated inclusion is denoted by:
$$\dd_{i}: \DDb^{n-1}\to \DDb^{n},\quad  i=0, \ldots ,n.$$
An $n+1$-simplex can be degenerated into an $n$-simplex by reducing the edge from $j$ to $j+1$ to a point. The associated surjection is denoted by:
$$\ss_{j}:\DDb^{n+1}\to \DDb^{n},\quad  j=0, \ldots ,n.$$

It is straightforward (and a good exercise) to check that these maps satisfy the following relations:
\begin{align*}
\dd_j\dd_i &= \dd_i\dd_{j-1},\quad i<j \\
\ss_j\ss_i &= \ss_i\ss_{j+1},\quad i\leq j \\
\ss_j\dd_i &= {\begin{cases} \dd_i\ss_{j-1}, & i<j \\
\id ,& i=j,\;i=j+1 \\ \dd_{i-1}\ss_j ,& i>j+1\end{cases}} \\
\end{align*}

It will prove helpful to think of the face maps $\dd_{i}$ and the degeneracy maps $\ss_{j}$ as morphisms in a category $\DD$, whose objects are the integers $[n]$ for $n=0, 1, \ldots $. Another way of presenting $\DD$ is as follows. Setting $[n]:=\{0, 1, \ldots, n\}$ (ordered finite set) a morphism in $\DD$ is a non-decreasing set map. The face map $\dd_{i}$ corresponds to the injective map which ``forgets'' $i$ in the range set. The degeneracy map $\ss_{j}$ corresponds to the surjective map which sends $j$ and $j+1$ to $j$  in the range set. Under this categorical interpretation we have a functor
$$\DDb^{\cdot}: \DD \to \texttt{Top}, \quad [n] \mapsto \DDb^{n},$$
where $\texttt{Top}$ is the category of topological spaces and continuous maps. A functor from $\DD$ to $\texttt{Top}$ is called a \emph{cosimplicial space}.

The opposite of the category $\DD$ is denoted by $\DD^{op}$, and the opposite of the face and degeneracy morphisms are denoted by:
$$d_{i}: [n] \to [n-1]\quad \textrm{and}\quad s_{j}: [n] \to [n+1], \quad i=0,\ldots , n.$$

They satisfy the \emph{simplicial relations}:
\begin{align*}
d_id_j &= d_{j-1}d_i,\quad i<j \\
s_is_j &= s_{j+1}s_i,\quad i\leq j \\
d_js_i &= {\begin{cases} s_{j-1}d_i ,& i<j \\
\id ,& i=j,\;i=j+1 \\ s_jd_{i-1} ,& i>j+1\end{cases}} \\
\end{align*}

A functor from $\DD^{op}$ to the category of sets, resp.\ modules, resp. topological spaces, is called a \emph{simplicial set}, resp.\ \emph{simplicial module}, resp.\ \emph{simplicial space}. It is sometimes helpful to think about a covariant functor from $\DD^{op}$ as a contravariant functor from $\DD$.

\subsection{Geometric realization}\label{geometricrealization} Given a simplicial set $X_{\cdot}: \DD^{op}\to \texttt{Set}$ one can form its \emph{geometric realization} as follows. First one takes a simplex of dimension $n$ for each element in the set $X_{n}$ (called the set of $n$-simplices), and then one glues them and contracts them according to the simplicial structure. More precisely one puts an equivalence relation on the union $\bigcup_{n} X_{n}\times \DDb^{n}$. This equivalence relation is 

$$(x,\varphi_{*}(y)) \sim (\varphi^{*}(x), y), \quad \textrm{for any}\quad  x\in X_{m}, y\in \DDb^{n},$$
and for any morphism $\varphi: [n] \to [m]$ in $\DD$.
The geometric realization of the simplicial set $X_{\cdot}$ is 
$$|X_{\cdot}|:= \bigcup_{n} X_{n}\times \DDb^{n}/ \sim .$$
 It gives rise to a well-defined functor from the category of simplicial sets to the category of topological spaces. One can easily extend this geometric realization to simplicial spaces by taking into account the topology of $X_{n}$ in the product $X_{n}\times \DDb^{n}$. Hence we obtain a functor from simplicial spaces to  topological spaces.

\subsection{Key example: the circle}\label{circle} Let us suppose that we want to realize the circle as the geometric realization of some simplicial set. The simplest way of constructing the circle as a cell complex is to start with an interval and to glue the two end points  together:
$$
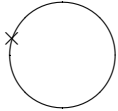
\xymatrix@=10pt{
&&&\\
&\, {\times}\ar@{-} `d/20pt[r]  `r[u] `u[l] `l[d]&&\\
&&&}$$

\bigskip

 Let us denote by $*$ this base-point, and by $\tau$ the $1$-cell.
The simplicial set $S^{1}_{\cdot}$ needs to have an element in $S^{1}_{0}$, that we denote by $*$, and an element in $S^{1}_{1}$, that we denote by $\tau$. The face maps on $\tau$ are obvious: $d_{0}(\tau)=d_{1}(\tau)= *$. But this is not enough because we want to build a simplicial set, and, in a simplicial set, we need to have degeneracies. So, in dimension $1$ we need to introduce the element $s_{0}(*)$. Similarly, in dimension $n$ we need to introduce the element $s_{0}{}^{n}(*)$ and the elements
$s_{n-1}\cdots \widehat{s_{i-1}} \cdots s_{0}(\tau)$ for $i= 1, \ldots , n$. This is enough because of the relations satisfied by the degeneracies.

So we see that not only $S^{1}_{n}$ has $n+1$ elements but we get naturally a unique bijection with the additive group $\ZZ/(n+1)\ZZ$:
$$ s_{0}{}^{n}(*) \leftrightarrow 0, \quad s_{n-1}\cdots \widehat{s_{i-1}} \cdots s_{0}(\tau) \leftrightarrow i\  \textrm{ for } i= 1, \ldots , n.$$
In low dimension we have

$$\begin{array}{l}
S^{1}_{3}=\{  *, s_{2}s_{1}\tau, s_{2}s_{0}\tau, s_{1}s_{0}\tau\}\\
S^{1}_{2}=\{  *, s_{1}\tau, s_{0}\tau\} \\
S^{1}_{1}=\{  *, \tau\} \\
S^{1}_{0}=\{  *\}  
\end{array}$$

\subsection{Example: classifying space of a discrete group} Let $G$ be a discrete group. One constructs a simplicial set $B_{\cdot}G$ as follows. In degree $n$ one puts $B_{n}G= G^{n}$ (cartesian product of $n$ copies of $G$) and the simplicial operations are given by

\begin{eqnarray*}
d_{0}(g_{1},\ldots, g_{n}) & = & (g_{2},\ldots, g_{n}),\\
d_{i}(g_{1},\ldots, g_{n}) & = & (g_{1},\ldots, g_{i}g_{i+1}, \ldots, g_{n})\quad \textrm{for } 1\leq i \leq n-1,\\
d_{n}(g_{1},\ldots, g_{n}) & = & (g_{1},\ldots, g_{n-1}),\\
s_{j}(g_{1},\ldots, g_{n}) & = & (g_{1},\ldots, g_{j}, 1, g_{j+1}, \ldots, g_{n})\quad \textrm{for } 0\leq j \leq n .
\end{eqnarray*}

It is a particular case of the nerve of a category, a notion that we do not need in this paper. The geometric realization of the simplicial set $B_{\cdot} G$ is a connected topological space denoted by $BG:= |B_{\cdot}G|$ and called the \emph{classifying space} of the discrete group $G$. It has the following homotopy properties: its fundamental group $\pi_{1}(BG)$ is the group $G$ and all the higher homotopy groups are trivial (this is because the discrete set $G$ has only $\pi_{0}$).

\subsection{Example: the free loop space}\label{cyclicbar}  Let $G$ be a discrete group. One constructs a simplicial set $\Gamma_{\cdot}G$, called the \emph{cyclic bar construction} of $G$,  as follows. In degree $n$ one puts $\Gamma_{n}G= G^{n+1}$ (cartesian product of $n+1$ copies of $G$) and the simplicial operations are given by

\begin{eqnarray*}
d_{i}(g_{0},\ldots, g_{n}) & = & (g_{0},\ldots, g_{i}g_{i+1}, \ldots, g_{n})\quad \textrm{for } 0\leq i \leq n-1,\\
d_{n}(g_{0},\ldots, g_{n}) & = & (g_{n}g_{0},\ldots, g_{n-1}),\\
s_{j}(g_{0},\ldots, g_{n}) & = & (g_{0},\ldots, g_{j}, 1, g_{j+1}, \ldots, g_{n})\quad \textrm{for } 0\leq j \leq n .
\end{eqnarray*}

\begin{thm}\label{cyclicbarloop} The geometric realization of the simplicial set $\Gamma_{\cdot} G$ is homotopy equivalent to the free loop space over the classifying space of $G$:
$$|\Gamma_{\cdot}G| \stackrel{\cong}{\longrightarrow} \mathcal{L} BG.$$
\end{thm}
Heuristically this result is not surprising since $\Gamma_{n}G= G^{n+1}=\Hom (S^{1}_{n}, G)$. But it is best understood in terms of cyclic sets and so postponed to the next section (see for example \cite{HC} Chapter 6). A straightforward generalization consists in starting with a simplicial group (resp.\ topological group) in place of a group. The standard geometric realization has to be replaced by the geometric realization of a bisimplicial set, resp.\ simplicial topological space.

\subsection{Geometric realization of a cosimplicial space}\label{cosimplicialrealization} Given a cosimplicial space $Z^{\cdot}$, whose dimension $n$ component is denoted by $Z(n)$ (in order to avoid confusion between cosimplicial dimension and iterated cartesian product), one constructs its geometric realization as being the subspace $$|| Z^{\cdot} || := \Hom_{\texttt{cosimpTop}}( \DDb^{\cdot}, Z^{\cdot})\quad \subset \quad \Pi_{n\geq 0} \Hom_{\texttt{Top}}(\DDb^{n}, Z(n)).$$

For any simplicial set $K_{\cdot}$ and any topological space $X$ one can form a cosimplicial space $\Hom(K_{\cdot}, X)$. There is a binatural map
$$|| \Hom(K_{\cdot}, X) || \to \Hom_{\texttt{Top}}(|K_{\cdot}|, X).$$

In good cases this map turns out to be a homeomorphism, see Theorem \ref{PatrasThomas}.

\subsection{Chain complex}\label{chaincplx} Let us start with a simplicial module $C_{\cdot}$ over some commutative ring $\KK$. Then we can form a chain complex $C_{\bullet}$ as follows. The module of $n$-chains is $C_{n}$ and the boundary map is 
$$b:= \sum_{i=0}^{n} (-1)^i d_{i}: C_{n}\to C_{n-1}.$$
The relation $b^{2}=0$ is a consequence of the relations between the face maps $d_{i}$. Observe that, at this point, the degeneracy maps are not used. They can be used to make this chain complex smaller. Indeed, let $C'_{n}$ be the submodule of $C_{n}$ spanned by the elements $s_{j}(x)$ for any $j$ and any $x\in C_{n-1}$. Then, the relations between the face maps and the degeneracy maps imply that $(C'_{\bullet}, b)$ is a subcomplex of $(C_{\bullet}, b)$ and is acyclic. Hence the quotient map 
$$C_{\bullet}\to \overline{C}_{\bullet} := C_{\bullet}/C'_{\bullet}$$
 is a quasi-isomorphism and we can use this quotient $ \overline{C}_{\bullet} $, called the \emph{normalized complex} of $C_{\cdot}$, to compute the homology of $C_{\bullet}$.

The tensor product of the simplicial modules $C_{\cdot}$ and $D_{\cdot}$ is a simplicial module $C_{\cdot}\t D_{\cdot}$ such that $(C_{\cdot}\t D_{\cdot})_{n}= C_{n}\t D_{n}$. However the tensor product of the complexes $C_{\bullet}$ and $D_{\bullet}$ is a chain complex $C_{\bullet}\t D_{\bullet}$ such that $(C_{\bullet}\t D_{\bullet})_{n}= \oplus_{p+q=n}C_{p}\t D_{q}$. The comparison is given by the so-called \emph{Alexander-Whitney map}:
$$
\xymatrix@R=4pt{
AW: &C_{n}\t D_{n}\ar[r] &  \bigoplus_{p+q=n}C_{p}\t D_{q}\\
&x\t y \ar@{|->}[r] & \sum_{p+q=n} {d_{max}}^{q}(x)\t  {d_{0}}^{p}(y)
}
$$
This chain map is known to be a quasi-isomorphism whose quasi-inverse is given by the so-called \emph{Eilenberg-Zilber map}, induced by the shuffles.

\subsection{Example: the singular complex}\label{singcplx} Let $X$ be a path-connected topological space. We let $S_{n}(X)= \Hom(\DDb^{n}, X)$ be the set of all continuous maps from the simplex to $X$. From the cosimplicial structure of $\DDb$ we get a simplicial structure on $S_{\cdot}(X)$. Taking the free $\KK$-module in each dimension we get a simplicial module $\mathcal{S}_{*}(X)$. By definition the homology of the associated chain complex (see above) is called the \emph{singular homology} of the topological space $X$.

The two functors ``geometric realization'' and ``singular simplicial set'' are adjoint to each other:
$$\xymatrix
{ \texttt{Top} \ar@<-1ex>[rr]_{S_{*}}&& \texttt{Simp Set}\ar@<-1ex>[ll]_{|-|}}
$$
 Moreover they have some good properties with respect to homology. In particular, if we start with a simplicial set $X_{\cdot}$, then we can make it into a simplicial module $\KK[X_{\cdot}]$ and then take the homology. On the other hand we can realize $X_{\cdot}$ to get a topological space $X=|X_{\cdot}|$ and take its singular homology. The unit of the adjunction induces an isomorphism (cf.\ for instance \cite{Curtis}):
$$ H_{\bullet}(\KK[X_{\cdot}]) \cong H_{\bullet}(\mathcal{S}_{*}(X)) =:H_{\bullet}(X).$$

\subsection{The singular cochain algebra} The diagonal map $\textrm{Diag}: X\to X\times X$ (we refrain to denote it by $\DD$ for obvious reasons) induces a coassociative product 
$$
\xymatrix{
\mathcal{S}_{*}(X)\ar[rr]^{\mathcal{S}_{*}(\textrm{Diag})}&&\mathcal{S}_{*}(X\times X)\ar[rr]^{AW} && \mathcal{S}_{*}(X)\t \mathcal{S}_{*}(X).
}$$
Since mathematicians better like to work with algebras rather than coalgebras, we dualize this object as follows. The \emph{cochains module} $\mathcal{S}^{n}(X):= \Hom(\mathcal{S}_{n}(X), \KK)$ is the linear dual of the chain module. The coalgebra structure on the chains, becomes an algebra structure on the cochains.  The resulting space  $\mathcal{S}^{*}(X)$ is a differential graded associative algebra (dga algebra for short). The cohomology of the underlying complex is, by definition, the \emph{singular cohomology} of the space $X$:
$$  H^{\bullet}(X):= H^{\bullet}(\mathcal{S}^{*}(X)).$$
The associative product on cochains induces a graded product on cohomology called   the \emph{cup-product} and denoted by $x\cup y$. It is known to be graded commutative, see for instance \cite{Greenberg}. 

\subsection{Example: Hochschild homology} Let $A$ be a unital associative algebra over $\KK$. We can form a simplicial module $C_{\cdot}(A)$ as follows:
$$ C_{n}(A):= A^{\t n+1}\ ,$$
$$\begin{array}{l}
d_{i}(a_{0}, \ldots , a_{n}) = (a_{0}, \ldots, a_{i}a_{i+1}, \ldots , a_{n}) \textrm{\ for\ } i=0, \ldots ,n-1,\\
d_{n}(a_{0}, \ldots , a_{n}) = (a_{n}a_{0}, a_{1}, \ldots, a_{n-1}),\\
s_{j}(a_{0}, \ldots , a_{n}) = (a_{0}, \ldots, a_{j}, 1, a_{i+1}, \ldots , a_{n}) \textrm{\ for\ } j=0, \ldots ,n.
\end{array}$$

The simplicial structure follows immediately from the associativity and unitality properties of $A$. By \ref{chaincplx} we get a chain complex, which is called the \emph{Hochschild complex} of $A$. Its homology is denoted by $HH_{\bullet}(A)$ (you will also find $HH_{\bullet}(A,A)$ and $H_{\bullet}(A,A)$ the literature, because we can replace the first occurence of $A$ in $A^{\t n+1}$ by any $A$-bimodule).

One can extend the Hochschild functors from the category of associative algebras to the category of \emph{differential graded associative algebras} (dga algebras for short). The trick is classical, for $A= \oplus_{m}A_{m}$ with boundary map $d$, the module $A^{\t n+1}$ is also graded and $d$ is extended by derivation to it. So we get a bicomplex with differentials induced by $b$ and $d$. The total complex is called the \emph{Hochschild complex of the dga algebra} $(A,d)$.

\subsection{The de Rham algebra} For any unital commutative algebra $A$, the module of \emph{K\" ahler differential forms} 
$\Omega^{1}_{A}$ is the quotient of $A\t A$ by the relation 
$$ ab\t c - a \t bc + ca \t b = 0 \ \textrm{ for any } a,b,c \in A.$$
It is customary to adopt Leibniz notation and denote the equivalence class of $a\t b$ by $a\, db$, so that 
$$d(ab) = a (db) + b (da).$$
We observe that $\Omega^{1}_{A}$ is an $A$-module and that $d1=0$. The module of higher differential forms is defined through the exterior power over $A$:
$$\Omega^{n}_{A}:= \Lambda^{n}_{A} (\Omega^{1}_{A}).$$
It is spanned by the elements $a_{0}da_{1}\cdots da_{n}, \textrm { for } a_{i}\in A$.

By definition the \emph{de Rham complex} of $A$ is the cochain complex
$$A= \Omega^{0}_{A}\to \Omega^{1}_{A}\to \cdots \to \Omega^{n}_{A} \stackrel{d}{\to} \Omega^{n+1}_{A} \to \cdots $$
where $d(a_{0}da_{1}\cdots da_{n})= da_{0}da_{1}\cdots da_{n}$. Its homology is denoted by $H^{n}_{DR}(A)$.

When $\KK=\RR$ and $A$ is the algebra of $C^{\infty}$-functions on a connected manifold $M$, then it is well-known that there is an isomorphism
$$H^{n}_{DR}(A)= H^{n}(M).$$
see for instance \cite{BottTu}.

The relationship between differential forms and Hochschild homology goes as follows. 

\begin{prop} If $\KK$ is a characteristic zero field, then the maps
$$ \epsilon : \Omega^{n}_{A} \to HH_{n}(A) ,\quad  a_{0}da_{1}\cdots da_{n}\mapsto \frac{1}{n!} \sum_{\ss\in S_{n}} \sgn (\ss) a_{0}da_{\ss(1)}\cdots da_{\ss(n)},$$
and
$$ \pi : HH_{n}(A) \to  \Omega^{n}_{A} ,\qquad (a_{0},a_{1},\cdots ,a_{n})\mapsto a_{0}da_{1}\cdots da_{n},$$
are well-defined and verify $\pi\circ \epsilon = \id$. Morevover, if $A$ is smooth, then they are both isomorphisms.
\end{prop}

This statement is classical and can be found in many places, including \cite{HC} Appendix E by Mar\' \i a Ronco. Let us recall that, here, $A$ being smooth means the following. For any commutative algebra $R$ equipped with an ideal $I$ such that $I^{2}=0$, we suppose that any algebra map $A \to R/I$ admits a lifting to $R$. The last assertion of the Proposition is known in the literature under the name ``Hochschild-Kostant-Rosenberg theorem''. In words it says that: homologically, smooth algebras behave like free algebras.

In the graded module $ \Omega^{*}_{A}:=\sum_{n\geq 0} \Omega^{n}_{A}$ we have $ \Omega^{0}_{A}=A$. The associative and commutative algebra structure of $A$ is extended into a graded commutative algebra structure on $ \Omega^{*}_{A}$ by
$$a_{0}da_{1}\cdots da_{p} * a'_{0}da'_{1}\cdots da'_{q} = a_{0}a'_{0}da_{1}\cdots da_{p}da'_{1}\cdots da'_{q}.$$
It can be shown that there is also a graded commutative algebra structure on the Hochschild homology of a commutative algebra compatible with the algebra structure. It is explicitly constructed out of the shuffles, cf.\ for instance \cite{HC} Chapter 3.

%%%%%%%%%%%%%%%%%%%%%%%%%%%%%%%%%%%%%%%%%%%%%%%%%%%%%%%%%%%%%%%%%%%%%%%%%%%%%%

\section{Hochschild homology and free loop space}

The aim of this section is to compare the Hochschild homology of the de Rham algebra of a simply connected manifold $M$ and the cohomology of its free loop space:
$$HH_{\bullet}(\Omega^{*}M)\cong H^{\bullet}(\mathcal{L} M)\ .$$
It is a consequence of an analogous statement where the de Rham algebra is replaced by the singular cochain algebra of a simply connected space $X$:

$$HH_{\bullet}(\mathcal{S}^{*}(X))\cong H^{\bullet}(\mathcal{L} X)\ .$$

The key object of the proof is the cocyclic space $\Hom (S^{1}_{\cdot}, X)$. Since $S^{1}_{n}= \ZZ/(n+1)\ZZ$ it gives rise to the Hochschild complex of the singular cochain algebra $\mathcal{S}^{*}(X)$ on one hand. On the other hand, taking the geometric realization gives the free loop space $\mathcal{L} X$.

We treat this case in the first part of this section. In the second part we treat a similar result involving the \emph{homology} of the free loop space instead of the cohomology.

\subsection{Homology of a cosimplicial space \cite{PatrasThomas}} For any  cosimplicial space  $Z^{\cdot}$ there are two ways to construct a dga algebra over $\KK$. One one hand one can take its geometric realization $|| Z^{\cdot} ||$ and then take the singular cochain algebra $\mathcal{S}^{*} || Z^{\cdot} || $. On the other hand the singular cochain algebras $\mathcal{S}^{*} Z^{q}
$ assemble into a cosimplicial-simplicial module (i.e.\ a functor $\DD\times \DD^{op}\to \texttt{Mod}$)  $\mathcal{S}^{*} Z^{*} $, giving rise to a bicomplex that we normalize in the simplicial direction to get $\overline{\mathcal{S}}^{\bullet} Z^{\bullet} $. By using the Eilenberg-Zilber maps (i.e.\ the shuffles) one can show that the total complex $\Tot\ \overline{\mathcal{S}}^{\bullet} Z^{\bullet} $ is a dga algebra. The comparison of these two dga algebras is given by a natural dga algebra map
$$\phi :\Tot\ \overline{\mathcal{S}}^{\bullet} Z^{\bullet}  \longrightarrow \mathcal{S}^{*}(|| Z^{\cdot} ||).$$

When $\phi$ happens to be a quasi-isomorphism, then the cosimplicial space $Z^{\cdot}$ is said to be \emph{convergent}. Let $K_{\cdot}$ be a simplicial set (for instance $S^{1}_{\cdot}$) and let $X$ be a topological space. Then $Z^{\cdot}:= \Hom( K_{\cdot}, X)$ is a cosimplicial space. The following convergence result is proved by F.\ Patras and J.-C.\ Thomas in \cite{PatrasThomas}  for $\KK$ being a field, enhancing the work of Bendersky and Gitler \cite{BenderskyGitler}.

\begin{thm}\label{PatrasThomas} Let $K_{\cdot}$ be a simplicial set whose geometric realization is a finite cell complex of dimension $\dim(K_{\cdot})$. Let $X$ be a topological space which is ${\textrm Conn}(X)$-connected (${\textrm Conn}(X)=0$ means connected, $\textrm{Conn}(X)=1$ means simply connected). If $\dim(K_{\cdot})\leq {\textrm Conn}(X)$, then $\Hom(K_{\cdot}, X)$ is a convergent cosimplicial space, that is:

$$\phi :\Tot\ \overline{\mathcal{S}}^{\bullet} \Hom(K_{\cdot}, X)  \longrightarrow \mathcal{S}^{*}(|| \Hom(K_{\cdot}, X)||)$$

is a quasi-isomorphism.
\end{thm}

\subsection{The cosimplicial space $\Hom (S^{1}_{\cdot}, X)$} Since $S^{1}_{\cdot}$ is a simplicial set, it follows that $\Hom (S^{1}_{\cdot}, X)$ is a cocyclic space. In dimension $n$ we get $X^{n+1}$ since $S^{1}_{n}= \ZZ/(n+1)\ZZ$. By \ref{cosimplicialrealization} the geometric realization of this cocyclic space is $$|| \Hom (S^{1}_{\cdot}, X) || = \Hom_{\texttt{Top}}(| S^{1}_{\cdot} |) , X)=: \mathcal{L} X .$$

The cosimplicial-simplicial module $\mathcal{S}^{*}(\Hom(S^{1}_{\cdot},X))$ is obtained by applying the functor $\mathcal{S}^{*}$ dimensionwise. The total complex $(\Tot\ M)_{\bullet}$ of a cosimplicial-simplicial module $M_{\cdot}({\cdot})$ (simplicial as index, cosimplicial in parenthesis) is obtained as follows:

$$(\Tot\ M)_{n}:= \bigoplus_{p-q=n} M_{p}(q),\quad  d = \sum_{i}(-1)^{i}d_{i} - \sum_{j}(-1)^{j}\dd_{j}$$
where $d_{i}: M_{p}(q)\to M_{p-1}(q)$ is a face of the simplicial structure and  $\dd_{j}: M_{p}(q)\to M_{p}(q+1)$ is a coface of the cosimplicial structure.

By the Alexander-Whitney map we obtain a morphism of cosimplicial-simplicial modules

$$ (\mathcal{S}^{*}(X))^{\t n+1} \longrightarrow \mathcal{S}^{*}(\Hom(S^{1}_{n},X)), n\geq 0,$$
which gives readily an isomorphism

$$ C_{\bullet}(\mathcal{S}^{*}(X)) \longrightarrow \Tot\ (\mathcal{S}^{*}(\Hom(S^{1}_{\cdot},X)) .$$

\begin{thm}\label{mainthm}  For any simply connected space $X$ there is a functorial isomorphism:

$$HH_{\bullet}(\mathcal{S}^{*}(X))\cong H^{\bullet}(\mathcal{L} X)\ .$$
\end{thm}
\begin{proof} We consider the composite

$$ C^{*}(\mathcal{S}^{*}(X)) \longrightarrow \Tot\ (\mathcal{S}^{*}(\Hom(S^{1}_{\cdot},X)) \longrightarrow \mathcal{S}^{*}(|| X ||).$$

Since $S^{1}= |S^{1}_{\cdot}|$ is a cellular complex of dimension $1$, by Theorem \ref{PatrasThomas} the map $\phi$ is a quasi-isomorphism as soon as $X$ is simply connected ($\dim S^{1} \leq \textrm{Conn}(X)$). Taking the homology of the complexes gives the result.  
\end{proof}

\subsection{Relationship with de Rham} If $M$ is a manifold, then its cohomology can be computed out of the de Rham algebra $\Omega^{*}(M)$. This de Rham algebra is a dga algebra which is quasi-isomorphic to the singular cochain algebra, see for instance \cite{BottTu}. By standard homological arguments, Hochschild homology of homotopy equivalent dga algebras are isomorphic. Hence we get $HH_{*}(\mathcal{S}^{*}(X))\cong HH_{*}(\Omega^{*}(M))$ and, as a corollary:

\begin{cor} For any simply connected manifold $M$ there is a functorial isomorphism:

$$HH_{\bullet}(\Omega^{*}(M))\cong H^{\bullet}(\mathcal{L} M)\ .$$

\end{cor}

\subsection{Homological version of Hochschild-free loop space relationship: cyclic bar construction}\label{HHandcyclicbar} There is also a nice relationship between the \emph{homology} (instead of cohomology) of the free loop  space with Hochschild homology. This time, the key object is the \emph{cyclic bar construction} $\Gamma_{\cdot}G$ of the discrete group $G$. From its definition given in  \ref{cyclicbar} it comes immediately that its associated simplicial module over $\KK$ is $[n] \mapsto A^{\t n+1}$ with $A=\KK[G]$, the group algebra of $G$ over $\KK$. Hence we get readily an isomorphism:

$$HH_{\bullet}(\KK[G]) \cong H_{\bullet}(\Gamma_{\cdot}G).$$

From the adjunction between the geometric realization of a simplicial set and the singular complex, we know that there is an isomorphism $H_{\bullet}(\Gamma_{\cdot}G) \cong H_{\bullet}(|\Gamma_{\cdot}G|)$. By Theorem \ref{cyclicbarloop} the space $|\Gamma_{\cdot}G|$ is homotopy equivalent to the free loop space of $BG$, whence an isomorphism
$$H_{\bullet}(|\Gamma_{\cdot}G|)\cong H_{\bullet}(\mathcal{L} BG).$$
Putting everything together we get the following

\begin{thm} For any discrete group $G$ there is a functorial isomorphism
$$HH_{\bullet}(\KK[G]) \cong  H_{\bullet}(\mathcal{L} BG).$$
\end{thm}

\subsection{Extension to simplicial groups}\label{extension}

Let $G_{\cdot}$ be a reduced simplicial group, i.e.\ $G_{n}$ is a group and $G_{0}={*}$. The cyclic bar construction can be extended to simplicial groups to give a bisimplicial set $\Gamma_{\cdot}G_{\cdot}$ (recall that a bisimplicial set is a functor $\DD^{op}\times \DD^{op} \to \texttt{Set}$). There are three ways to realize a bisimplicial set:
\begin{itemize}
 \item by realizing the first simplicial structure, getting a simplicial set, then realizing it,
 \item by realizing the second simplicial structure, getting a simplicial set, then realizing it,
 \item by realizing the diagonal of the bisimplicial set (which is a simplicial set).
 \end{itemize}
 By \cite{Quillen66} we know that the three methods end up with spaces which are homotopy equivalent.
 
 On the chain complex side, a bisimplicial module gives rise to a chain bicomplex. Its homology is, by definition, the homology of the total complex of this bicomplex.
 
 The aforementioned results extend to the simplicial framework straightforwardly and we get 
\begin{thm} For any simplicial  group $G_{\cdot}$ there is a functorial isomorphism
$$HH_{\bullet}(\KK[G_{\cdot}]) \cong  H_{\bullet}(\mathcal{L} BG_{\cdot}).$$
\end{thm}

The interest of this generalization lies in the fact that any (reasonable) connected topological space  (e.g. a  connected manifold) is homotopy equivalent to the geometric realization of a reduced simplicial group.

%%%%%%%%%%%%%%%%%%%%%%%%%%%%%%%%%%%%%%%%%%%%%%%%%%%%%%%%%%%%%%%%%%%%%%%%%%%%%%%%%%

\section{Cyclic}\label{s:cyclic}

In this section we recall the notion of cyclic objects and their properties. See for instance \cite{LodayQuillen, HC} for details. 

\subsection{Cyclicity of the simplicial circle}\label{cyclicity} In \ref{circle} we have seen that the set of $n$-simplices of the circle can be identified with the cyclic group of order $n+1$. Let us adopt a multiplicative notation for this group, with generator $\tau_n:= s_{n-1}\ldots s_{1}(\tau) \in S^1_n$ and relation $\tau_n{}^{n+1}= \id_n := s_0{}^{n}(*)$. In particular $\tau_0=*$ and $\tau_1=\tau$. From this definition it is easy to find the relations between these cyclic operators and the face and degeneracy operations. We get:

$$\begin{array}{l}
\tau_n \dd_i = \dd_{i-1} \tau_{n-1}, \quad \textrm{for } 1\leq i\leq n,\\
\tau_n \dd_0 = \dd_{n} \\
\tau_n \ss_i = \ss_{i-1}  \tau_{n+1}, \quad \textrm{for } 1\leq i\leq n,\\
\tau_n \ss_0 = \ss_{n} \tau_{n+1}{}^2 .
\end{array}
$$

We are now ready to introduce a new category, denoted by $\DD C$ and called the \emph{cyclic category} or \emph{Connes' category}. It is generated by the face operations, the degeneracy operations and the cyclic operations, subject to all the relations mentioned above. It has the following property.

\begin{prop} In the category $\DD C$, whose objects are the sets $[n]$, $n\geq 0$,  any morphism can be uniquely written as $\phi \circ g$ where $\phi$ is a morphism of $\DD$ and $g$ is an element of the cyclic group (i.e. $\Aut ([n])\cong \ZZ/(n+1)\ZZ$).
\end{prop}

Observe that the Proposition justifies the notation.
By definition a \emph{cocyclic set}, resp.\ \emph{cocyclic space}, is a functor $\DD C\to \texttt{Set}$, resp.\ $\DD C\to \texttt{Top}$.

\subsection{Cyclic objects} The opposite category of $\DD C$ is the category $\DD C ^{op}$ generated by the morphisms:
$$\begin{array}{l}
d_i:[n]\to [n-1],  \quad \textrm{for } 1\leq i\leq n,\\
s_j:[n]\to [n+1],  \quad \textrm{for } 1\leq i\leq n,\\
t_n : [n] \to [n] ,
\end{array}
$$
satisfying the relations 

\begin{align*}
d_id_j &= d_{j-1}d_i,\quad i<j \\
s_is_j &= s_{j+1}s_i,\quad i\leq j \\
d_js_i &= {\begin{cases} s_{j-1}d_i, & i<j \\
\id , & i=j,\;i=j+1 \\ 
s_jd_{i-1}, & i>j+1\end{cases}} \\
d_i t_n &= t_{n-1}d_{i-1} , \quad  1\leq i\leq n\\
d_0 t_n &= d_{n} ,\\
s_i t_n &=   t_{n+1}s_{i-1}, \quad 1\leq i\leq n\\
s_0 t_n &=  t_{n+1}{}^2s_{n}\\
(t_n)^{n+1}&= \id_n
\end{align*}

By definition a \emph{cyclic set}, resp.\ \emph{cyclic space}, resp.\ \emph{cyclic module}, is a functor from the category $\DD C ^{op}$ to the category  $\texttt{Set}$, resp.\ $\texttt{Top}$, resp.\ \texttt{Mod}.

\subsection{On the $S^{1}$-structure of the geometric realization of a cyclic set} It is well-known that the circle can be equipped with a group structure once it is identified with $SO(2)$. If we see it as the interval $[0,1]$ with the two endpoints identified, then the group structure is given by the addition of real numbers modulo 1. We will show that this group structure can be read off from the cyclic structure of $S^{1}_{\cdot}$. In fact, for any cyclic set $X_{\cdot}$ the geometric realization $|X_{\cdot}|$ is equipped with an action of the topological group $S^{1}$. It would have been helpful if this structure were induced by a simplicial (or even cyclic) map $S^{1}_{\cdot}\times X_{\cdot}\to X_{\cdot}$. This is not the case. What holds is the following. From the cyclic structure of $X_{\cdot}$ we will be able to construct an intermediate cyclic set $F(X)_{\cdot}$ together with simplicial maps
$$S^{1}_{\cdot}\times X_{\cdot} \leftarrow F(X)_{\cdot} \to X_{\cdot}\quad .$$
When passing to the geometric realization the left one becomes a homeomorphism, hence one can take its inverse and we get the expected action of $S^{1}$ on $X= |X_{\cdot}|$:
$$ S^{1} \times X \cong |F(X)_{\cdot}| \to X.$$

This intermediate space is constructed as follows. There is an obvious forgetful functor from cyclic sets to simplicial sets. The functor 
$$F: \texttt{Simp Set} \to  \texttt{Cyclic Set}$$
 is simply its left adjoint. Explicitly it is constructed as follows. If $Y_{\cdot}$ is a simplicial set, then $F(Y)_{n}= S^{1}_{n}\times Y_{n}$ and the cyclic structure is given by:
 
 $$
 \begin{array}{l}
 f_{*}(g,y) =(f_{*}(g), (g^{*}(f))_{*}(y) \\
 h^{}{*}(g,y)= ( hg, y)\\
 \textrm{ for any } f \textrm{ in } \DD^{op}, \ g \textrm{ and } h \in S^{1}_{n}\quad .
 \end{array}
 $$

So, in plain words, it is essentially the product of $S^{1}_{\cdot}$ by $Y_{\cdot}$, but where the cyclic structure has been twisted. If $X_{\cdot}$ is a cyclic set, then there exists an obvious map (evaluation)
$$ ev: F(X)_{\cdot}\to X_{\cdot}$$
which consists in evaluating the cyclic group action. The two forgetful maps $F(X)_{\cdot}\to S^{1}_{\cdot}$ and $F(X)_{\cdot}\to X_{\cdot}$ are obviously simplicial morphisms. It is proved in details in \cite{HC} Chapter 7 that the resulting map $|F(X)_{\cdot}| \to S^{1} \times X$ is a homeomorphism and that the composite $S^{1}\times X\to X$ does give an action of the topological group $S^{1}$ on $X$. This result has been first noted by several authors independently: Burghelea and Fiedorowicz \cite{BurgheleaFiedo}, Goodwillie \cite{Goodwillie}, J.D.S. Jones \cite{Jones87} for instance.  

\subsection{Examples of cyclic objects} The examples of simplicial objects given in section 1 are also cyclic objects.

\noindent$\bullet$ \emph{Classifying space of a discrete group}.  Let $G$ be a discrete group and let $z\in G$ be an element in the center of $G$. Then the map
$$t_{n}(g_{1}, \ldots , g_{n}) := (z(g_{1}g_{2}\cdots g_{n})^{-1}, g_{1}, \ldots, g_{n-1}) $$
makes $B_{\cdot}G$ into a cyclic set.

\noindent{\sc Exercise.} It is known that the classifying space of the discrete group $\ZZ$, that is the geometric realization of $B_{\cdot}\ZZ$, is homotopy equivalent to $S^{1}$. Show that $B_{\cdot}\ZZ$ is a cyclic set and provide a cyclic (hence simplicial) map $S^{1}_{\cdot}\to B_{\cdot}\ZZ$ realizing this homotopy equivalence.

\noindent$\bullet$ \emph{The cocyclic space $\Hom(S^{1}_{\cdot},X)$}. Since $S^{1}_{\cdot}$ is a cyclic set the cosimplicial space $\Hom(S^{1}_{\cdot},X)$ is a cocyclic space. The resulting action of the topological group $S^{1}$ on the free loop space $\mathcal{L} X=|| \Hom(S^{1}_{\cdot},X) ||$ is by rotating loops.

\noindent$\bullet$ \emph{The cyclic bar construction.} The simplicial set $\Gamma_{\cdot}G$ defined in \ref{cyclicbar} is a cyclic set for the cyclic action given by
$$t_{n}(g_{0}, \ldots , g_{n}) := (g_{n},  g_{0}, \ldots , g_{n-1}).$$
The resulting action of   $S^{1}$ on the free loop space $\mathcal{L} BG$ is the expected one (rotating loops).

\noindent$\bullet$ \emph{ The cyclic module of an algebra}.  Let $A$ be an associative algebra. On $A^{\t n+1}$ the action of $t_{n}$ is given by
$$t_{n}(a_{0}, \ldots, a_{n}):= (a_{n}, a_{0}, \ldots, a_{n-1}).$$
Sometimes it is technically helpful to modify this action by multiplying by the signature, that is $(-1)^{n}$.

\subsection{Cyclic modules and cyclic homology theories} First, let us recall that for the (multiplicative) cyclic group ${\bf C}_{n}=\{ t=t_{n} | t_{n}{}^{n+1} = 1\}$ there is a periodic free resolution of the trivial module:

$$\KK   \longleftarrow \KK[{\bf C}_{n}]   \stackrel{1-t}{\longleftarrow} \KK[{\bf C}_{n}]  \stackrel{N}{\longleftarrow} \KK[{\bf C}_{n}]   \stackrel{1-t}{\longleftarrow} \KK[{\bf C}_{n}]  \stackrel{N}{\longleftarrow} \cdots $$
where $N:= 1 + t + \cdots + t^{n}$ is the norm map.

The good relationship between the cyclic operator and the face operators in a cyclic module like $[n] \mapsto A^{\t n+1}$ enables us to build a bicomplex $CC_{**}(A)$ called the \emph{cyclic bicomplex}:

$$\xymatrix{
*{}\ar[d] & *{}\ar[d] & *{}\ar[d] & *{}\ar[d] & \\
A^{\t 3}\ar[d]_{b} & A^{\t 3}\ar[l]_{1-t}\ar[d]_{-b'} & A^{\t 3}\ar[l]_{N}\ar[d]_{b} & A^{\t 3}\ar[l]_{1-t}\ar[d]_{-b'} & \ar[l]_{N} \\
A^{\t 2}\ar[d]_{b} & A^{\t 2}\ar[l]_{1-t}\ar[d]_{-b'} & A^{\t 2}\ar[l]_{N}\ar[d]_{b} & A^{\t 2}\ar[l]_{1-t}\ar[d]_{-b'} & \ar[l]_{N} \\
A& A\ar[l]_{1-t}& A\ar[l]_{N}& A\ar[l]_{1-t}& \ar[l]_{N} \\
}
$$

\medskip

\noindent where $b:= \sum_{i=0}^{n} (-1)^i d_{i}$ and $b':= \sum_{i=0}^{n-1} (-1)^i d_{i}$. In order to make it work we have to modify slightly the action of $t_{n}$ in the linear case by multiplying it by the signature, i.e.\ by $(-1)^{n}$, see \cite{HC} for details.

When $A$ is unital, every other columns (the one with $b'$) are acyclic because there exists an explicit homotopy $(a_{0},a_{1},\cdots ,a_{n})\mapsto (1, a_{0},a_{1},\cdots ,a_{n})$. The other columns are copies of the Hochschild complex. A careful analysis of the spectral sequence of the cyclic bicomplex leads to the following \emph{Connes periodicity exact sequence}:
$$\cdots \to HH_{n}(A) \stackrel{I}{\to} HC_{n}(A) \stackrel{S}{\to} HC_{n-2}(A) \stackrel{B}{\to} HH_{n-1}(A) \stackrel{I}{\to} \cdots $$

\subsection{Negative and periodic cyclic homology} From the periodic resolution of a cyclic group, any module $M$ over it gives  rise to the following chain complexes:
$$
\xymatrix{
 & & & M\ar@{<-}[r]^{1-t}  & M\ar@{<-}[r]^{N}  &M\ar@{<-}[r]^{1-t}  & \cdots\\
  \cdots \ar@{<-}[r] &M\ar@{<-}[r]^{1-t}  & M\ar@{<-}[r]^{N}  &M\ar@{<-}[r]^{1-t}  &M\ar@{<-}[r]^{N} &M\ar@{<-}[r]^{1-t}    &  \cdots\\
   \cdots \ar@{<-}[r] &M\ar@{<-}[r]^{1-t}  & M\ar@{<-}[r]^{N}  &M &&&\\
}$$

The common copy of $M$ is in degree $0$. Similarly the cyclic bicomplex $CC(A)$ of an algebra $A$ has two siblings denoted by $CC^{per}(A)$ and $CC^{-}(A)$, giving rise to the \emph{periodic cyclic homology} $HC^{per}(A)$ (sometimes denoted by $\widehat{HC}(A)$ by influence of Tate) and the \emph{negative cyclic homology} theory $HC^{-}(A)$. It is clear from the definitions of $HH, HC, HC^{per}, HC^{-}$ that they are related by various exact sequences, see \cite{HC} for details.

\subsection{The Gysin sequence of an $S^{1}$-space} The classifying space of the circle is the base-space of a well-known fibration:

$$S^{1} \to ES^{1} \to BS^{1}$$
(apply the functor $B$ to the fibration $\ZZ\to \RR \to S^{1}$, so $ES^{1}=B\RR$). 
More generally any $S^{1}$-space $X$ is the fiber of a fibration
$$
X \to ES^{1}\times_{S^{1}}X \to BS^{1}.
$$
The (co)homology of the total space is, by definition, the \emph{equivariant $S^{1}$-(co)homology} of the space $X$, denoted by $H_{n}^{S^{1}}(X)$ (resp.\ $H_{S^{1}}^{n}(X)$). The homology spectral sequence of this fibration degenerates into an exact sequence called the \emph{Gysin sequence}:

$$\cdots \to H_{n}(X) \to H_{n}^{S^{1}}(X)\to H_{n-2}^{S^{1}}(X)\to H_{n-1}(X)\to \cdots $$

There is a similar statement in cohomology.

To any $S^{1}$ space $Z$ (e.g. $\mathcal{L} X$)  one can associate not only the cohomology $H^{\bullet}(Z)$ and  the equivariant cohomology $H^{\bullet}_{S^{1}}(Z)$, but also
$\widehat{H}^{\bullet}_{S^{1}}(Z)$ and $G^{\bullet}_{S^{1}}(Z)$, which mimick the periodic (or Tate) cyclic homology and the negative cyclic homology respectively. 

%%%%%%%%%%%%%%%%%%%%%%%%%%%%%%%%%%%%%%%%%%%%%
\section{Cyclic homology and the free loop space} 

In this section we mix the results of sections 2 and 3 to compare the Hochschild-cyclic homology theories with the homology theories of the free loop space. 

\subsection{Cyclic version of the comparison map} We have seen in section 2 that there is an isomorphism  of graded modules $HH_{\bullet}(\SS^{* }X) \cong H^{\bullet}(\LL X)$ induced by the comparison of cosimplicial-simplicial modules. The precise isomorphism in degree $n$ is $HH_{-n}(\SS^{* }X) \cong H^{n}(\LL X)$ when we choose the ``homological version'' of the $\ZZ$-graded total complex associated to $C_{\bullet}(\SS^{* }X)$.
It turns out that these modules are \emph{cyclic} modules (see section 3), and that the comparison map is a cyclic map. So we can replace simplicial by cyclic. Hence it is not surprizing that a cyclic version of Theorem \ref{mainthm} gives:

\begin{thm}[\cite{Jones87}]  For any simply connected space $X$ there is a functorial isomorphism:

$$HC_{\bullet}^{-}(\mathcal{S}^{*}(X))\cong H^{\bullet}_{S^{1}}(\mathcal{L} X)\ .$$

\end{thm}

The reason for the appearance of $HC^{-}$, where we could have expected $HC$, comes from the use of cosimplicial objects in place of simplicial objects.

\begin{thm}[\cite{Jones87}]  For any simply connected space $X$ there is a commutative diagram:

\xymatrix@C=6pt{
\cdots \ar[r] & HH_{-n}(\mathcal{S}^{*}(X))\ar[r]\ar[d] &HC_{-n}^{-}(\mathcal{S}^{*}(X))\ar[r]\ar[d] & HC_{-n-2}^{-}(\mathcal{S}^{*}(X))\ar[r]\ar[d] &HH_{-n-1}(\mathcal{S}^{*}(X))\ar[r] \ar[d]& \cdots \\
\cdots \ar[r] &  H^{n}(\mathcal{L} X)\ar[r] & H^{n}_{S^{1}}(\mathcal{L} X)\ar[r] &  H^{n+2}_{S^{1}}(\mathcal{L} X)\ar[r] & H^{n+1}(\mathcal{L} X)\ar[r] & \cdots \\
}
\end{thm}

 There are, of course, similar results in the homological framework. They can be found in \cite{HC}.

\subsection{About the action of $O(2)$} So far we have only considered the action of the topological group $S^{1}=SO(2)$, but there is in fact an action of $O(2)$ on the free loop space. The discretization of $O(2)$ involves the dihedral groups in place of the cyclic groups. The point is the following: there is a periodic resolution for the cyclic groups (of period 2), but there is no periodic resolution for the dihedral groups. However there is a periodic resolution for the \emph{quaternionic groups} (of period 4), from which we can take advantage of to work out the $O(2)$-action. We refer to \cite{Loday87, Lodder90} and \cite{HC} Chapter 5 for the details, and to \cite{Spalinski10} for the most recent article on this theme.

\subsection{Other algebraic properties}  In section 3 we have exploited the fact that the simplicial objects used in section 2 are, in fact, cyclic objects. But there are more algebraic properties at hand. For instance the cohomology of a space is a graded commutative algebra. So it is natural to ask oneself how this algebraic structure on $H^{*}(\mathcal{L} X)$ is reflected on $HH_{*}(\mathcal{S}^{*}(X))$. This question has been solved in \cite{Menichi01, Menichi09} by L.\ Menichi.

In the seminal paper \cite{ChasSullivan} M.\ Chas and D.\ Sullivan introduced a product
$$H_{p}(\mathcal{L}(M)) \t H_{q}(\mathcal{L}(M)) \to H_{p+q-n}(\mathcal{L}(M)) $$
for any manifold $M$ of dimension $n$. This discovery prompted a new theory, called ``String Topology'' which let appear a lot of algebraic structures like: Batalin-Vilkovisky algebras, their homotopy version \cite{DCV}, the cacti operad \cite{CohenVoronov}, involutive Lie bialgebras, and many more. In \cite{CohenJones} R.\ Cohen and J.\ Jones study this \emph{loop product} on the Hochschild homology under the isomorphism of section 2. 

One can find in \cite{CohenVoronov} self-contained notes treating the string topology operations   from different perspectives, such as geometric intersection theory, operads and Morse theory.

%%%%%%%%%%%%%%%%%%%%%%%%%%%%%%%%%%%%%%%%%%%%

\section{Conclusion and comments} In order to understand the homological and homotopical properties of the free loop space of a manifold, we have used simplicial techniques. First, it should be said, that, if one is interested only in the rational homotopy and simply connected framework, then working with the Sullivan minimal model is good enough, since one knows how to construct the minimal model of $\mathcal{L} M$ out of the minimal model of $M$, cf. \cite{SullivanVigue}.

On top of the simplicial techniques it is clear that the bar-cobar adjoint functors form a key tool. But we have only used them in the framework of associative (co)algebras, that is algebras and coalgebras over the nonsymmetric operad $As$ (see for instance \cite{LodayValletteAO} Chapter 2). However, since we want to take into account the circle action, we need to see the associative (co)algebras over $As$, but where $As$ is viewed as a ``cyclic operad''. Let us recall here that when one wants to handle the bar-cobar constructions for commutative algebras, one needs to work also with Lie algebras (think of rational homotopy theory), since the Koszul dual of the operad $Com$ is the operad $Lie$. In this case we have the advantage that the Koszul dual of an operad is still an operad. Now the questions that we face here are the following. First,  what is the Koszul dual of a cyclic operad ? Second, what is the Koszul dual of $As$ considered as a cyclic operad ? Third, play the bar-cobar game. 

The next step will be to take into account, in these computations, the ``symmetry up to homotopy'' properties.

\end{document}